\documentclass[11pt]{article}

\usepackage{amsmath}
\usepackage{amsfonts}

\newsymbol\Subset 1362

\newcommand{\be}{\begin}
\newcommand{\Ref}[1]{(\ref{#1})}

\newcommand{\cb}{\mathcal B}
\def\cc{\mathcal C}

\newcommand{\cg}{\mathcal G}

\def\cl{\mathcal L}

\def\fl{\mathcal R}

\newcommand{\al}{\alpha}
\newcommand{\del}{\delta}

\newcommand{\eps}{\epsilon}
\def\ga{\gamma}

\def\ka{\kappa}

\def\La{\Lambda}
\def\om{\omega}
\def\Om{\Omega}
\newcommand{\si}{\sigma}

\newcommand{\co}{\mathbb C}

\newcommand{\e}{{\mathbb E}}

\def\bl{\mathbb L}
\def\bs{\mathbb S}
\def\bt{\mathbb T}

\newcommand{\R}{\mathbb R}

\newcommand{\la}{\langle}
\newcommand{\ra}{\rangle}
\newcommand{\st}{\,:\,}
\newcommand{\ti}{\tilde}

\newcommand{\prtl}{\partial}
\def\square{\kern20pt{\vbox{\hrule height.4pt
        \hbox{\vrule width.4pt height 6pt\kern6pt
                \vrule width.4pt}
        \hrule height.4pt}}}

\newcommand{\rpy}{\R_+\times Y}

\def\inv{^{-1}}
\def\hol{\text{Hol}}

\def\rf{{\R}^4}
\def\proof{{\em Proof.}\ }

\def\SW{\text{SW}}

\def\U#1{\text{U}(#1)}

\def\spc{\text{spin}^c}

\def\ind{\text{ind}}

\def\hx{{\mathfrak b}}

\def\ft{\mathfrak t}

\def\sW{\mathsf W}

\def\llw#1#2#3{{L^{#1,#2}_#3}}
\def\lw#1#2{{L^{#1,#2}}}

\def\endt#1{_{:#1}}

\def\tu{^{(T)}}
\def\xt{{X\tu}}
\def\mt{M\tu}
\def\gt{G\tu}

\def\rr{_o}

\def\hatf{\hat F}

\newcommand{\tcg}{\breve\cg}

\def\zd{Z'}
\def\oxk{\Om^+_{X,K}}
\def\emm{\vphantom{I}_\eps M}
\def\Em#1{\vphantom{I}^#1\!M}
\def\emb{\emm_\hx}
\def\emxy{\emm_{x_0,x_1}}
\def\onm{\vphantom{I}_1M}
\def\omred{\omega_{\text{red}}}
\def\ared{{A_{\text{red}}}}
\def\son{_{x_1}}
\def\sz{_{x_0}}
\def\szo{_{x_0,x_1}}

\begin{document}

\title{A generalized blow-up formula for\\Seiberg--Witten invariants}

\author{Kim A.\ Fr\o yshov}

\date{19 April 2006}

\maketitle

\bibliographystyle{plain}

\be{abstract}
We prove a gluing formula for Seiberg--Witten invariants
which describes in particular
the behaviour of the invariant under blow-up and rational blow-down.
\end{abstract}


\section{Introduction}
\label{sec:intro}

In this paper we use the results of \cite{Fr10,Fr11} to establish
a gluing formula for Seiberg--Witten invariants of certain $4$--manifolds
containing a negative definite piece. The formula describes in particular
the behaviour of the Seiberg--Witten invariant under blow-up and under
the rational blow-down procedure introduced by Fintushel--Stern~\cite{FS9}.
While the formula will be known in principle
to experts, to our knowledge
no complete proof has previously been published. As for the 
classical blow-up formula, this was proved by
Bauer~\cite[Corollary~4.2]{Bauer2}. (Earlier, a proof had been announced
by Salamon~\cite{das2}, and there is a sketch of a proof in
Nicolaescu~\cite{nico1}.)
In the case of rational blow-down the formula 
was stated by Fintushel--Stern~\cite[Theorem~8.5]{FS9} with
a brief outline of a proof.
Apart from providing a proof in the general case,
the main motivation for writing this paper was to show how the parametrized 
version of the gluing theorem in \cite{Fr11}
can be used to handle at least the simplest cases
of obstructed gluing, thereby providing a unified approach to a wide range
of gluing problems.

Before stating our results we explain how the
Seiberg--Witten invariant, usually defined for closed $4$--manifolds,
can easily be generalized to compact $\spc$ $4$--manifolds
$Z$ whose boundary $Y'=\prtl Z$ satisfies $b_1(Y')=0$ and admits a metric $g$
of positive scalar curvature. (By a $\spc$--manifold we mean as in \cite{Fr10}
an oriented smooth manifold with a $\spc$--structure.)
As usual we also assume that
$b_2^+(Z)>1$. Let $\{Y'_j\}$ be the components of $Y'$, which are
rational homology spheres.
Let $\ti Z$ be the manifold with tubular ends
obtained from $Z$ by adding a half-infinite tube $\rpy'$.
Choose a Riemannian metric on $\ti Z$
which agrees with $1\times g$ on the ends. We consider the monopole equations
on $\ti Z$ perturbed solely by means of a smooth $2$--form $\mu$ on $\ti Z$
supported in $Z$ as in \cite[Equation~13]{Fr10}.
Let $M=M(\ti Z)$ denote the moduli space of monopoles
over $\ti Z$ that are asymptotic over $\rpy'_j$ to the unique (reducible)
monopole over $Y'_j$.
For generic $\mu$ the moduli space $M$
will be free of reducibles and a smooth compact manifold of dimension
\[\dim M=2h(Y')+\frac14(c_1(\cl_{Z})^2-\si(Z))-b_0(Z)+b_1(Z)-b_2^+(Z),\]
see \cite[Section~9]{Fr10}.
Choose a base-point $x\in\ti Z$ and let $M_x$
be the framed moduli space defined just as $M$ except that
we now only divide out by those gauge transformations $u$ for which
$u(x)=1$. Let $\bl\to M$ be the complex line bundle whose
sections are given by maps $s:M_x\to\co$ satisfying
\be{equation}\label{eqn:U1equiv}
s(u(\om))=u(x)\cdot s(\om)
\end{equation}
for all $\om\in M_x$ and gauge transformations $u$.
A choice of homology 
orientation of $Z$ determines an orientation of $M$, and
we can then define the Seiberg--Witten invariant of $Z$ just as for closed
$4$--manifolds:
\[\SW(Z)=\be{cases}
\la c_1(\bl)^k,[M]\ra & \text{if}\quad\dim M=2k\ge0,\\
0 & \text{if}\quad\text{$\dim M$ is negative or odd}.
\end{cases}\]
The use of $\bl$ rather than $\bl\inv$ prevents a sign in
Theorem~\ref{thm:gen-blow-up} below. (Another justification is that, although
$M_x\to M$ is a principal bundle with respect to the canonical $\U1$--action,
it seems more natural to regard that action as a {\em left} action.)
This invariant $\SW(Z)$ depends only on the homology oriented
$\spc$--manifold $Z$, not on the choice of
positive scalar curvature metric $g$ on $Y'$; the proof of this is a special
case of the proof of the generalized blow-up formula, which we are now ready
to state.

\be{thm}\label{thm:gen-blow-up}
Let $Z$ be a connected, compact, homology oriented
$\spc$ $4$--manifold
whose boundary $Y'=\prtl Z$ satisfies $b_1(Y')=0$ and admits a metric of
positive scalar curvature, and such that $b_2^+(Z)>1$.
Suppose $Z$ is separated by an embedded rational homology sphere $Y$ admitting
a metric of positive scalar curvature,
\[Z=Z_0\cup_YZ_1,\]
where $b_1(Z_0)=b_2^+(Z_0)=0$. Let $Z_1$ have the orientation, homology
orientation, and $\spc$ structure inherited from $Z$. Then
\[\SW(Z)=\SW(Z_1)\quad\text{if}\quad \dim M(\ti Z)\ge0.\]
\end{thm}

We will show in Section~\ref{sec:prelim} that $\dim M(\ti Z_0)\le-1$.
(A particular case of this was proved by different methods in 
\cite[Lemma~8.3]{FS9}.)
The addition formula for the index then yields
\[\dim M(\ti Z)=\dim M(\ti Z_0)+1+\dim M(\ti Z_1)\le\dim M(\ti Z_1).\]

The following corollary describes the effect on the Seiberg--Witten invariant
of both ordinary blow-up and rational blow-down:

\be{cor}
Let $Z_0,\zd_0,Z_1$ be compact, homology oriented
$\spc$ $4$--manifolds with $-\prtl Z_1=\prtl Z_0=\prtl\zd_0=Y$
as $\spc$ manifolds, where $Y$ is a $\spc$ rational
homology sphere admitting a metric of positive scalar curvature. Suppose
$b^+_2(Z_1)>1$, $b_1(Z_0)=b_1(\zd_0)=0$, and $b_2(Z_0)=b_2^+(\zd_0)=0$. Let
\[Z=Z_0\cup_YZ_1, \quad \zd=\zd_0\cup_YZ_1\]
have the orientation, homology orientation and $\spc$ structure induced
from $Z_0,\zd_0,Z_1$. Then
\[\SW(Z)=\SW(\zd)\quad\text{if}\quad\dim M(\zd)\ge0.\]
\end{cor}

\proof Set $n_\pm=\dim M(\pm\ti Z_0)$ and $W=Z_0\cup_Y(-Z_0)$. Then
\[-1=\dim M(W)=n_++1+n_-,\]
hence $n_\pm=-1$. Thus
\[\dim M(Z)=\dim M(\ti Z_1)\ge \dim M(\zd)\ge0.\]
The theorem now yields
\[\SW(Z)=\SW(Z_1)=\SW(\zd).\square\]

\section{Preliminaries on negative definite $4$--manifolds}
\label{sec:prelim}

Let $X$ be a connected $\spc$ Riemannian $4$--manifold
with tubular ends
$\rpy_j$, $j=1,\dots,r$, as in \cite[Subsection~1.3]{Fr10}. Suppose
each $Y_j$ is a rational homology sphere and $b_1(X)=0=b_2^+(X)$.
We consider the monopole equations on $X$ perturbed only by means of a
$2$--form $\mu$ as in \cite[Equation~13]{Fr10}, where now
$\mu$ is supported in a given non-empty, compact, codimension~$0$ submanifold
$K\subset X$. Let $\al_j\in\fl_{Y_j}$ be the reducible monopole over $Y_j$
and $M_\mu=M(X;\vec\al;\mu;0)$ the moduli space of monopoles over $X$
with asymptotic limits $\vec\al=(\al_1,\dots,\al_r)$. This moduli space
contains a unique reducible point $\om(\mu)=[A(\mu),0]$. Let $\oxk$ 
denote the space of (smooth) self-dual $2$--forms on $X$ supported in $K$,
with the $C^\infty$ topology. Let $p$ and $w$ be the exponent and
weight function used in the definition of the moduli space $M_\mu$, as
in \cite[Subsection~3.4]{Fr10}.

\be{lemma}\label{lemma:ind-neg-pos}
Let $R$ be the set of all $\mu\in\oxk$ such that the operator
\be{equation}\label{eqn:dirac-mu}
D_{A(\mu)}:\llw pw1(\bs^+_X)\to\lw pw(\bs^-_X)
\end{equation}
is either injective or surjective. Then $R$ is open and dense in $\oxk$.
\end{lemma}

Of course, whether the operator is injective or surjective for a given
$\mu\in R$ is determined by its index, which is independent of $\mu$.

\proof By \cite[Proposition~2.2~(ii)]{Fr10} and the proof of
\cite[Proposition~5.2]{Fr10}, the operator
\[d^+:\ker(d^*)\cap\llw pw1\to\lw pw\]
is an isomorphism. Therefore, if $A\rr$ is a reference
connection over $X$ with limits $\al_j$
as in \cite[Subsection~3.4]{Fr10} then there is a unique (smooth)
$a=a(\mu)\in\llw pw1$ with 
\[d^*a=0,\qquad d^+a=-\hatf^+(A\rr)-i\mu.\]
Hence we can take $A(\mu)=A\rr+a(\mu)$. Since the operator \Ref{eqn:dirac-mu}
has closed image, it follows by continuity of the map $\mu\mapsto A(\mu)$
that $R$ is open in $\oxk$.

To see that $R$ is dense, fix $\mu\in\oxk$ and write $A=A(\mu)$. Let $W$ be
a Banach space of smooth $1$--forms on $X$ supported in $K$ as provided by
\cite[Lemma~8.2]{Fr10}. Using the unique continuation property of the Dirac
operator it is easy to see that $0$ is a regular value of the smooth map
\be{align*}
h:W\times(\llw pw1(\bs^+_X)\setminus\{0\})&\to\lw pw(\bs^-_X),\\
(\eta,\Phi)&\mapsto D_{A+i\eta}\Phi.
\end{align*}
In general, if $f_1:E\to F_1$ and $f_2:E\to F_2$ are surjective homomorphisms
between vector spaces then $f_1|_{\ker f_2}$ and $f_2|_{\ker f_1}$ have
identical kernels and isomorphic cokernels. In particular, the projection
$\pi:h\inv(0)\to W$ is a Fredholm map whose index at every point agrees with
the index $m$ of $D_A$. By the Sard-Smale theorem the regular values of $\pi$
form a residual (hense dense) subset of $W$. If $\eta\in W$ is a regular value
then we see that $D_{A+i\eta}$ is injective when $m\le0$ and surjective when
$m>0$. Since the topology on $W$ is stronger than the $C^\infty$ topology
it follows that $R$ contains points of the form $\mu+d^+\eta$
arbitrarily close to $\mu$.\square

\be{lemma}\label{lemma:empty-inj}
Suppose the metric on each $Y_j$ has positive scalar curvature.
Let $R'$ be the set of all $\mu\in\oxk$ such that the irreducible part
$M^*_\mu$ is empty
and the operator $D_{A(\mu)}$ in \Ref{eqn:dirac-mu} is injective.
Then $R'$ is open and dense in $\oxk$.
\end{lemma}

\proof Recall that $M^*_\mu$ has expected dimension $2m-1$, where 
$m=\ind_{\co}D_{A(\mu)}$.

Suppose $m>0$. We will show that this leads to a contradiction.
Let $R''$ be the set of all $\mu\in\oxk$ for which $M_\mu$ is regular.
(Note that the reducible point is regular precisely when $D_{A(\mu)}$ is
surjective.)
From Lemma~\ref{lemma:ind-neg-pos} and \cite[Proposition~8.2]{Fr10} one finds
that $R''$ is dense in $\oxk$. (Starting with a given $\mu$, first perturb
it a little to make the reducible point regular, then a little more to
make also the irreducible part regular.)
But for any $\mu\in R''$ the moduli space
$M_\mu$ would be compact with one reducible point, which yields a
contradiction as in \cite{Fr1}.
Therefore, $m\le0$.

We now see, exactly as for $R''$, that $R'$ is dense in $\oxk$.
To prove that $R'$ is open we use a compactness argument
together with the following fact:
For any given $\mu_0\in R'$ there is a neighbourhood $U$ of
$\om(\mu_0)$ in $\cb(X;\vec\al)$ such that
\[M^*_\mu\cap U=\emptyset\]
for any $\mu\in\oxk$ with $\|\mu-\mu_0\|_p$ sufficiently small.
To prove this we work in a slice at $(A(\mu_0),0)$, ie we represent
$\om(\mu)$ (uniquely) by $(A,0)$ where $d^*(A-A(\mu_0))=0$, and we consider
a point in $M^*_\mu$ represented by $(A+a,\phi)$ where $d^*a=0$.
Note that since $b_1(X)=0$, the latter representative is unique up to
multiplication of $\phi$ by unimodular constants.

Observe that there is a constant $C_1<\infty$ such that if
$\|\mu-\mu_0\|_p$ is sufficiently small then
\[\|\psi\|_{\llw pw1}\le C_1\|D_A\psi\|_{\lw pw}\]
for all $\psi\in\llw pw1$. Hence if $L=(d^*+d^+,D_A)$ then for such $\mu$
one has
\[\|s\|_{\llw pw1}\le C_2\|Ls\|_{\lw pw}\]
for all $s\in\llw pw1$. Denoting by $\SW_\mu$ the Seiberg--Witten map
over $X$ for the perturbation form $\mu$ we have
\[0=\SW_\mu(A+a,\phi)-\SW_\mu(A,0)=(d^+a-Q(\phi),D_A\phi+a\phi),\]
where $Q$ is as in \cite{Fr10}. Taking $s=(a,\phi)$ we obtain
\[\|s\|_{\llw pw1}\le C_2\|Ls\|_{\lw pw}\le
C_3\|s\|^2_{\lw{2p}w}\le C_4\|s\|_{\llw pw1}^2.\]
Since $s\neq0$ we  conclude that
\[\|s\|_{\llw pw1}\ge C_4\inv.\]
Choose $\del\in(0,C_4\inv)$ and define
\[U=\{[A(\mu_0)+b,\psi]\st \|(b,\psi)\|_{\llw pw1}<\del,\:d^*b=0\}.\]
If $\|\mu-\mu_0\|_p$ is so small that
$\|A-A(\mu_0)\|_{\llw pw1}\le C_4\inv-\del$ then
\[\|(A+a-A(\mu_0),\phi)\|_{\llw pw1}
\ge\|s\|_{\llw pw1}-\|A-A(\mu_0)\|_{\llw pw1}\ge\del,\]
hence $[A+a,\phi]\not\in U$.\square 

\section{The extended monopole equations}
\label{sec:ext-mon-eqn}

We now return to the situation in Theorem~\ref{thm:gen-blow-up}.
Set $X_j=\ti Z_j$ for $j=0,1$. Choose metrics of positive scalar curvature
on $Y$ and $Y'$ and a metric on the disjoint union $X=X_0\cup X_1$ which 
agrees with the corresponding product metrics on the ends.
Let $Y$ be oriented
as the boundary of $Z_0$, so that $X_0$ has an end $\rpy$ and $X_1$
an end $\R_+\times(-Y)$. Gluing these two ends of $X$
we obtain as in \cite[Subsection~1.4]{Fr10} a manifold $\xt$ for
each $T>0$.

Choose smooth monopoles $\al$ over $Y$ and $\al'_j$ over $Y'_j$ (these
are reducible, and unique up to gauge equivalence). Let $S\rr=(A\rr,\Phi\rr)$
be a reference configuration over $X$ with these limits over the ends,
and $S'\rr$ the associated reference configuration over $\xt$.
Adopting the notation introduced in the beginning of
\cite[Subsection~2.2]{Fr11},
let $\cc$ be the corresponding $\llw pw1$ configuration space over $X$ and
$\cc'$ the corresponding $\llw p\ka1$ configuration space over $\xt$.
For any finite subset $\hx\subset X\endt0=Z_0\cup Z_1$ let
$\cg_\hx,\cg'_\hx$
be the corresponding groups of gauge transformations that restrict to $1$
on $\hx$.

As in Section~\ref{sec:prelim} we first consider the monopole equations
over $X$ and $\xt$ perturbed only by means of a self-dual $2$--form 
$\mu=\mu_0+\mu_1$, where $\mu_j$ is supported in $Z_j$. The corresponding
moduli spaces will be denoted $M(X)$ and $M\tu=M(\xt)$. Of course, 
$M(X)$ is a product of moduli spaces over $X_0$ and $X_1$:
\[M(X)=M(X_0)\times M(X_1).\]
By Lemma~\ref{lemma:empty-inj} we can choose $\mu_0$ such that $M(X_0)$
consists only of the reducible point (which we denote
by $\omred=[\ared,0]$), and such that the operator
\be{equation}\label{eqn:dared}
D_\ared:\llw pw1(\bs^+_{X_0})\to\lw pw(\bs^-_{X_0})
\end{equation}
is injective.
By \cite[Proposition~8.2]{Fr10} and unique continuation for self-dual closed
$2$--forms we can then choose $\mu_1$ such that
\be{itemize}
\item $M(X_1)$ is regular and contains no reducibles, and
\item the irreducible part of $\mt$ is regular for all natural numbers $T$.
\end{itemize}
Set
\[k=-\ind_{\co}(D_\ared)\ge0.\]
If $k>0$ then $\omred$ is not a regular point of $M(X_0)$
and we cannot appeal to the gluing theorem \cite[Theorem~2.1]{Fr11}
for describing $\mt$ when $T$ is large. We will therefore introduce
an extra parameter $z\in\co^k$ into the Dirac equation on $Z_0$, to obtain
what we will call the ``extended monopole equations'', such that $\omred$
becomes a regular point of the resulting parametrized
moduli space over $X_0$.
This will allow us to apply the gluing theorem for parametrized
moduli spaces, \cite[Theorem~5.1]{Fr11}.

We are going to add to the Dirac equation an extra term $\beta(A,\Phi,z)$ which
will be a product of three factors:
\be{description}
\item[(i)]a holonomy term $h_A$ (to achieve gauge equivariance)
\item[(ii)]a cut-off function $g(A,\Phi)$ (to retain an apriori pointwise
bound on $\Phi$)
\item[(iii)]a linear combination $\sum z_j\psi_j$ of certain negative spinors
(to make $\omred$ regular).
\end{description}
We will now describe these terms more precisely.

{\bf (i)} Choose an embedding
$f:\rf\to\text{int}(Z_0)$, and set $x_0=f(0)$ and $U_0=f(\rf)$. For each
$x\in U_0$ let $\ga_x:[0,1]\to U_0$ be the path from $x_0$ to $x$ given by
\[\ga_x(t)=f(tf\inv(x)).\]
For any
$\spc$ connection $A$ over $U_0$ define the function $h_A:U_0\to\U1$ by
\[h_A(x)=\exp\left(-\int_{[0,1]}\ga_x^*(A-\ared)\right),\]
cf.\ \cite[Equation~1]{Fr11}. Note that $h_A$ depends on the choice of
$\ared$, which is only determined up to modification by elements of $\cg$.

{\bf (ii)} Set $K_0=f(D^4)$, where $D^4\subset\R^4$ is the closed unit disk.
Choose a smooth function $g:\cb^*(K_0)\to[0,1]$ such that $g(A,\Phi)=0$ when
$\|\Phi\|_{L^\infty(K_0)}\ge2$
and $g(A,\Phi)=1$ when $\|\Phi\|_{L^\infty(K_0)}\le1$. Extend $g$ to 
$\cb(K_0)$ by setting $g(A,0)=1$ for all $A$.

{\bf (iii)} By unique continuation for the formal adjoint $D^*_\ared$ there are
smooth sections $\psi_1,\dots,\psi_k$ of $\bs^-_{X_0}$ supported in $K_0$
and spanning a linear complement of the image of the operator $D_\ared$
in \Ref{eqn:dared}.

For any configuration $(A,\Phi)$ over $X$ and $z=(z_1,\dots,z_k)\in\co^k$
define
\[\beta(A,\Phi,z)=g(A,\Phi)\,h_A\sum_{j=1}^kz_j\psi_j.\]
Note that for gauge transformations $u$ over $X$ one has
\[u(x_0)\,h_{u(A)}=u\,h_A.\]
Since $g$ is gauge invariant, this yields
\[\beta(u(A),u\,\Phi,u(x_0)z)=u\cdot\beta(A,\Phi,z).\]
The following lemma is useful for estimating the holonomy term $h_A$:

\be{lemma}\label{lemma:ha-est}
Let $a=\sum a_jdx_j$ be a $1$--form on the closed unit disk $D^n$ in $\R^n$,
$n>1$.
For each $x\in D^n$ let $J(x)$ denote the integral of $a$ along the line
segment from $0$ to $x$, ie
\[J(x)=\sum_{j=1}^nx_j\int_0^1a_j(tx)\,dt.\]
Then for any $q\ge1$ and $r>qn$
and non-negative integer $k$ there is a constant $C<\infty$
independent of $a$ such that
\[\|J\|_{L^q_k(D^n)}\le C\|a\|_{L^r_k(D^n)}.\]
\end{lemma}

\proof If $b$ is a function on $D^n$ and $\chi$ the characteristic
function of the interval $[0,1]$ then
\[\be{aligned}
\int_{D^n}\int_0^1b(tx)\,dt\,dx
&=\int_{D^n}b(x)\int_0^1t^{-n}\chi(t\inv|x|)\,dt\,dx\\
&=\frac1{n-1}\int_{D^n}(|x|^{1-n}-1)\,b(x)\,dx.
\end{aligned}\]
From this basic calculation the lemma is easily deduced.\square

It follows from the lemma that $a\mapsto h_{\ared+a}$ defines a smooth map
$L^p_1(K_0;i\La^1)\to L^q_1(K_0)$
provided $p>4q>16$. Hence, if $p>16$ (which we henceforth assume) then
\[\cc(K_0)\times\co^k\to L^p(K_0;\bs^-),
\quad((A,\Phi),z)\mapsto\beta(A,\Phi,z)\]
is a smooth map whose derivative at every point is a compact operator.
Here $\cc(K_0)$ is the $L^p_1$ configuration space over $K_0$.

The extended monopole equations for $((A,\Phi),z)\in\cc\times\co^k$ are
\be{equation}\label{eqn:ext-mon}
\be{gathered}
\hatf^+_A+i\mu-Q(\Phi)=0,\\
D_A\Phi+\beta(A,\Phi,z)=0.
\end{gathered}
\end{equation}
(Cf.\ the holonomy perturbations of the instanton equations constructed
in \cite[2\,(b)]{D2}.) We define 
actions of $\cg$ and $\cg'$ on $\cc\times\co^k$ and $\cc'\times\co^k$
respectively by
\[u(S,z)=(u(S),u(x_0)z).\]
Then the left hand side of \Ref{eqn:ext-mon} describes a $\cg$--equivariant
smooth map $\cc\times\co^k\to\lw pw$.

For $\eps>0$ let $B^{2k}_\eps\subset\co^k$ denote the open
ball of radius $\eps$ about the origin,
and $D^{2k}_\eps$ the corresponding closed ball.
For $0<\eps\le1$ set
\[\emb(X)=\left.\{\text{solutions $((A,\Phi),z)\in\cc\times B^{2k}_\eps$ to
\Ref{eqn:ext-mon}}\}\right/\cg_\hx,\]
This moduli space is clearly a product of moduli spaces over $X_0$ and
$X_1$:
\[\emb(X)=\emm_{\hx_0}(X_0)\times M_{\hx_1}(X_1),\]
where $\hx_j=\hx\cap X_j$.

Noting that the equations \Ref{eqn:ext-mon} also make sense over $\xt$ we
define
\[\emb\tu=\left.\{\text{solutions $((A,\Phi),z)\in\cc'\times B^{2k}_\eps$ to
\Ref{eqn:ext-mon}}\}\right/\cg'_\hx.\]

We define $\Em\eps_\hx(X)$ and $\Em\eps\tu_\hx$ in a similar way as
$\emb(X)$ and $\emb\tu$, but with $D^{2k}_\eps$ in place of $B^{2k}_\eps$.

Choose a base-point $x_1\in Z_1$. We will only consider the cases when
$\hx$ is a subset of $\{x_0,x_1\}$, and we indicate $\hx$ by listing
its elements (writing $\emxy$ and $\emm$ etc).



\be{lemma}Any element of $\Em1(X_0)$ or $\Em1\tu$ has a smooth representative.
\end{lemma}

\proof Given Lemma~\ref{lemma:ha-est} this is proved in the usual way.\square

\be{lemma}\label{lemma:phi-bound}
There is a $C<\infty$ independent of $T$ such that $\|\Phi\|_\infty<C$
for all elements $[A,\Phi,z]$ of $\Em1(X)$ or
$\Em1\tu$.
\end{lemma}

\proof Suppose $|\Phi|$ achieves a local maximum $\ge2$ at some point
$x$. If $x\not\in K_0$ then one obtains a bound on $|\Phi(x)|$ using the
maximum principle as in \cite[Lemma~2]{KM4}.
If $x\in K_0$ then the same works because then $g(A,\Phi)=0$.\square

\be{lemma}
$\Em1(X)$ and $\Em1\tu$ are compact for all $T>0$.
\end{lemma}

\proof Given Lemmas~\ref{lemma:ha-est} and \ref{lemma:phi-bound},
the second approach to compactness
in \cite{Fr10} carries over.\square

We identify $M_{\hx_0}(X_0)$ with the set of elements
of $\onm_{\hx_0}(X_0)$ with $z=0$, and similarly for moduli spaces over
$X,\xt$. It is clear from the definition of $\beta(A,\Phi,z)$ that $\omred$
is a regular point of $\onm(X_0)$. Since $\onm_{x_0}(X_0)$ has expected 
dimension $0$,
it follows that $\omred$ is an isolated point of $\onm_{x_0}(X_0)$. Because
$\Em1_{x_0}(X_0)$ is compact, there is an $\eps$ such that
$\emm_{x_0}(X_0)$ consists only of the point $\omred$. Fix such an $\eps$
for the remainder of the paper.

\be{lemma}\label{lemma:emm-ch-conv}
If $\om_n\in\emm^{(T_n)}$ with $T_n\to\infty$ then a subsequence of
$\{\om_n\}$ chain-converges to $(\omred,\om)$ for some $\om\in M(X_1)$.
\end{lemma}

\proof Again, this is proved as in \cite{Fr10} using the second approach
to compactness.\square


\be{cor}\label{cor:no-red}
If $T\gg0$ then $\emm\tu$ contains no element which is reducible
over $Z_1$.\square
\end{cor}

\section{Applying the gluing theorem}


Let $\hol=\hol_1$ be defined as in
\cite[Equation~1]{Fr11} in terms of a path in $\xt$
from $x_0$ to $x_1$ running once through the neck.

By \cite[Proposition~2.3]{Fr11}, if $K_1=(X_1)\endt\ft$ with $\ft\gg0$ then
there is a $\U1$--invariant open subset
$V_1\subset\cb^*\son(K_1)=\cb\son(K_1)$ containing
$R_{K_1}(M\son(X_1))$, and
a $\U1$--equivariant smooth map
\[q_1:V_1\to M\son(X_1)\]
such that
$q_1(\om|_{K_1})=\om$ for all $\om\in M\son(X_1)$.
Here $R_{K_1}$ denotes restriction to $K_1$.
It follows from Lemma~\ref{lemma:emm-ch-conv} that if $T$ is sufficiently large
then  $\om|_{K_1}\in V_1$ for all $\om\in\emm\tu\son$.

\be{prop}\label{prop:appl-gl-thm}
For all sufficently large $T$ the moduli space
$\emm\tu\son$ is regular and the map
\be{equation}\label{eqn:emmtuson}
\emm\tu\son\to M\son(X_1),\quad\om\mapsto q_1(\om|_{K_1})
\end{equation}
is an orientation preserving $\U1$--equivariant diffeomorphism.
\end{prop}

\proof We will apply the version of \cite[Theorem~5.1]{Fr11}
with (in the notation of \cite{Fr11})
$\bt$ acting non-trivially on $\sW$. Set
\be{gather*}
G=\emm\szo(X)=\{\omred\}\times M\son(X_1),\\
K=K_0\cup K_1,\\
V=\cb\sz(K_0)\times V_1\times B^{2k}_\eps.
\end{gather*}
Note that $G$ is compact and $\tcg_\hx(K)=\cg_\hx(K)$.
Define 
\[q:V\to G,\quad(\om_0,\om_1,z)\mapsto(\omred,q_1(\om_1)).\]
In general, an element $(u_0,u_1)\in\U1^2$ acts on appropriate
configuration and
moduli spaces like any gauge transformation $u$ with $u(x_j)=u_j$, $j=0,1$,
and it acts on $B^{2k}_\eps$ by multiplication with $u_0$.
Then clearly, $q$ is $\U1^2$--equivariant, so
by the gluing theorem there is a compact, codimension~$0$ submanifold
$K'\subset X$ containing $K$ and a $\U1^2$--equivariant open subset
$V'\subset\cb^*\szo(K')\times B^{2k}_\eps$ containing $R_{K'}(G)$ 
and satisfying $R_K(V')\subset V$ and such that
for all sufficiently large $T$ the space 
\[\gt=\{(\om,z)\in\emm\tu\szo\st(\om|_{K'},z)\in V'\}\]
consists only of regular points, and the map 
\be{equation}\label{eqn:emm-map}
\gt\to\U1\times\emm\szo(X),\quad(\om,z)
\mapsto(\hol(\om),(\omred,q_1(\om|_{K_1})))
\end{equation}
is a $\U1^2$--equivariant diffeomorphism.
But it follows from
Lemma~\ref{lemma:emm-ch-conv} that $\gt=\emm\tu\szo$ for $T\gg0$, and
dividing out by the action of $\U1\times\{1\}$ in \Ref{eqn:emm-map}
we see that \Ref{eqn:emmtuson} is a $\U1$--equivariant diffeomorphism.

We now discuss orientations. Given $\del=\pm1$ we will say a map is
{\em $\del$--preserving} if it changes 
orientations by the factor $\del$. Set
\[c:=b_1(X)+b^+_2(X).\]
By \cite[Proposition~4.3~(ii) and Theorem~5.1]{Fr11} the map \Ref{eqn:emm-map}
is $(-1)^{c+1}$--preserving. 
Using \cite[Proposition~4.4]{Fr11} it is a simple exercise to show that
$\emm\szo(X)\to M\son(X_1)$ is $(-1)^c$--preserving.
On the other hand, $(u_0,1)\in\U1\times\{1\}$ acts
on $\U1$ in \Ref{eqn:emm-map} by multiplication with $u_0\inv$.
Thus we have got three signs, which cancel each other since $(c+1)+c+1$
is even. Therefore,
the map \Ref{eqn:emmtuson} does preserve orientations.\square

{\em Proof of Theorem~\ref{thm:gen-blow-up}:} For large $T$ let
$\bl\to\emm\tu$ be the complex line bundle associated to the base-point
$x_1$ as in Section~\ref{sec:intro}. For $j=1,\dots,k$ the map 
\[s_j:\emm\tu\son\to\co,\quad[A,\Phi,z]\mapsto\hol(A)\cdot z_j\]
is $\U1$--equivariant in the sense of \Ref{eqn:U1equiv} and 
therefore defines a smooth section of $\bl$. The sections
$s_j$ together form a section $s$ of the bundle
$\e=\oplus^k\bl$ whose zero set
is the unparametrized moduli space $\mt$.
It is easy to see that $s$ is a regular section precisely when $\mt$
is a regular moduli space, which by Corollary~\ref{cor:no-red} and the
choice of $\mu_1$ holds at least when
$T$ is a sufficiently large natural number. In that case $s\inv(0)=\mt$
as oriented manifolds.  Set 
\[\ell=\frac12\dim\,\mt\ge0,\]
so that $\dim\,M(X_1)=2(k+\ell)$. If $\ell$ is not integral then
$\SW(Z_1)=0=\SW(Z)$ and we are done. Now suppose $\ell$ is integral and
let $T$ be a large natural number.
Choose a smooth section $s'$ of $\e'=\oplus^\ell\bl$ such that
$\si=s'|_{\mt}$ is a regular section of $\e'|_{\mt}$, or equivalently,
such that $s\oplus s'$ is a regular section of
$\e\oplus\e'=\oplus^{k+\ell}\bl$. Then
\[\SW(Z_1)=\#(s\oplus s')\inv(0)=\#\si\inv(0)=\SW(Z),\]
where the first equality follows from Proposition~\ref{prop:appl-gl-thm},
and $\#$ as usual means a signed count.
\square

\noindent\textsc{Fakult\"at f\"ur Mathematik, Universit\"at Bielefeld,\\
Postfach 100131, D-33501 Bielefeld, Germany.}\\
\\
E-mail:\ froyshov@math.uni-bielefeld.de

\end{document}